\documentstyle{amsart}

\numberwithin{equation}{section}

\newtheorem{theorem}{Theorem}[section]
\newtheorem{theorema}{Theorem}
\newtheorem{lemma}[theorem]{Lemma}

\newtheorem{proposition}[theorem]{Proposition}

\theoremstyle{remark}
\newtheorem{remark}[theorem]{Remark}

\newtheorem{ack}{Acknowledgment}

\begin{document}

\title{Ricci-flat K\"ahler metrics on canonical bundles}
\author{Roger Bielawski}
\thanks{Research supported by an EPSRC Advanced Research Fellowship}


\address{Department of Mathematics\\
University of Glasgow\\
Glasgow G12 8QW\\
Scotland }

\email{rb@@maths.gla.ac.uk}

\begin{abstract} We prove the existence of a (unique) $S^1$-invariant Ricci-flat K\"ahler metric on a neighbourhood of the zero section in the canonical bundle of a  real-analytic K\"ahler manifold $X$, extending the metric on $X$.    \end{abstract}

\maketitle

In the important paper \cite{Cal}, Calabi proved existence of Ricci-flat K\"ahler metrics on two classes of manifolds: a) cotangent bundles of projective spaces; b) canonical bundles of K\"ahler-Einstein manifolds. The metrics on $T^\ast {\Bbb C}P^n$ are actually hyperk\"ahler and in the intervening years hyperk\"ahler metrics were shown to exist on cotangent bundles of many other K\"ahler manifolds. Finally, recently, B. Feix \cite{Feix} and, independently, D. Kaledin \cite{Kha} have shown that a real-analytic K\"ahler metric on a complex manifold $X$ always extends to a (essentially unique) hyperk\"ahler metric on a neighbourhood of $X$ in $T^\ast X$.
\par
The aim of this paper is to prove the analogous generalization for the other class of Calabi's metrics. Our main existence result can be stated as follows:

\begin{theorema} Let $X$ be a real-analytic K\"ahler manifold. Then there exists a unique Ricci-flat K\"ahler metric on a neighbourhood of $X$ in the canonical bundle $K_X$ of $X$ which extends the metric on $X$ and for which the standard $S^1$-action on $K_X$ is isometric and Hamiltonian. \label{main}\end{theorema}
 The condition of real-analycity of the K\"ahler metric is clearly necessary, since the extended metric is Ricci-flat.\\
We also notice that the adjunction formula shows that the canonical bundle is the only line bundle over $X$ which can admit a Ricci-flat K\"ahler metric.

\section{Proof of Theorem 1}

Let $M$ be an $n+1$ dimensional K\"ahler manifold with a free Hamiltonian circle action. Then the metric can be locally 
written in the form:
\begin{equation}G=\sum g_{ij}dz_i\otimes d\bar{z}_j +wdt^2 +w^{-1}\phi^2, \label{metric}\end{equation}
where $t$ is the moment map on $M$, $\phi$ is the circle-invariant 1-form and the $z_i$ are local coordinates on 
$M/{\Bbb C}^\ast$.\\
The complex structure $I$ maps $dt$ to $w^{-1}\phi$. Pedersen and Poon \cite{PP} (and LeBrun \cite{LeB} for $n=1$) have worked out the conditions for the complex structure to be integrable and for the metric to be Einstein (in fact, Pedersen and Poon deal with the more general case of torus symmetry). We recall their theorem.
\begin{theorem}{\bf [Pedersen-Poon]} Let $w$ be a smooth positive funtion and $[g_{ij}]$ a positive definite hermitian matrix of smooth functions on an open set $U$ in ${\Bbb C}^n\times {\Bbb R}$. The metric \eqref{metric} is Ricci-flat if and only if the following system of equations holds for some constant $c$:
\begin{equation} 4u_{z_i\bar{z}_j}+c(g_{ij})_t=0,\label{PP}\end{equation}
\begin{equation} u_t=cw. \label{cw}\end{equation}
\begin{equation} 4w_{z_i\bar{z}_j}+(g_{ij})_{tt}=0\label{new}\end{equation}
Here $u$ is defined by 
\begin{equation} \det g=we^u.\label{u}\end{equation}
\par
Furthermore, the metric is defined on a circle bundle over $U$ if and only if the cohomology class $[F]$ of the curvature of $\phi$ which is given by
\begin{equation} F=-\left(\frac{i}{2}(g_{ij})_tdz_i\wedge d\bar{z}_j+iw_{z_i}dt\wedge dz_i-iw_{\bar{z}_j}dt\wedge d\bar{z}_j\right) \label{curvature}\end{equation}
belongs to $2\pi {\Bbb Z}$.\hfill $\Box$
\label{P-P}\end{theorem}

The constant $c$ has the following significance:
\begin{proposition} $\Delta_G t= c$.\end{proposition}
\begin{pf} For any function $f$ we have
\begin{equation} \Delta_G f= g^{ij}  \left( 4\frac{\partial^2f}{\partial z_i\partial \bar{z}_j}+
w^{-1}\frac{\partial f}{\partial t} \frac{\partial g_{ij}}{\partial t}\right) +\frac{\partial}{\partial t}\left( \frac{\partial f}{\partial t}w^{-1}\right). \label{Laplacian}\end{equation}

Thus, for $f=t$, we obtain
$$\Delta_G t= g^{ij}w^{-1}\frac{\partial g_{ij}}{\partial t} + \frac{\partial w^{-1}}{\partial t}= g^{ij}\frac{\partial w^{-1} g_{ij}}{\partial t}= w^{-1} \frac{\partial\ln\det g}{\partial t}+ \frac{\partial w^{-1}}{\partial t}. $$
Now, using \eqref{u} and \eqref{cw}, we have
$$\Delta_G t= w^{-1} \frac{\partial\ln\det g}{\partial t}+ \frac{\partial w^{-1}}{\partial t} =w^{-1}\frac{\partial\ln w}{\partial t} +c + \frac{\partial w^{-1}}{\partial t}=c.$$ 
\end{pf}

For hyperk\"ahler manifolds, the constant $c$ can take only two values:

\begin{proposition} Let $M^{4n}$ be a hyperk\"ahler manifold equipped with an isometric  and Hamiltonian (for one symplectic structure) action of the circle. Then the moment map $t$ for this action is harmonic if the action is tri-holomorphic and satisfies $\Delta t=n$ otherwise. \end{proposition}
\begin{pf} If the action is triholomorphic, then the moment map is the real part of a complex moment map. If the action is not triholomorphic, i.e. it rotates the complex structures orthogonal to $I$, then the moment map is a K\"ahler potential for another K\"ahler form (corresponding to the complex structure $J$) and the result follows. \end{pf}

 We shall seek metrics with $c\neq 0$. In this case the equation \eqref{new} is the consequence of the other two equations. Moreover we can eliminate the function $w$ from the equations and replace \eqref{cw} and \eqref{u} with
\begin{equation} \left(e^u\right)_t=c\det g.\label{toda2}\end{equation}
\medskip

Suppose now that we are given a K\"ahler metric $h=\sum h_{ij}dz_i\otimes d\bar{z}_j$ on a complex $n$-dimensional manifold $X$ and we wish to extend $h$ to a Ricci-flat K\"ahler metric $g$ in a neighbourhood of $X$ in a line bundle $L$. Furthermore we require that the canonical $S^1$-action on $L$ is Hamiltonian. Clearly, a necessary condition for this is that we can solve the following (singular) Cauchy problem:
\begin{equation}\left\{\begin{matrix} u_{x_i{x}_j}+u_{y_i{y}_j}+c(g_{ij})_t=0\\
  (e^u)_t =c\det g\quad\quad \quad\quad\\ \left(g_{ij}\right)_{|_{t=0}} =h_{ij} \quad\quad\quad\quad\\
  \left(e^u\right)_{|_{t=0}} = 0\quad\quad\quad\quad .\end{matrix}\right.\label{system}\end{equation}
Here $z_i=x_i+\sqrt{-1}y_i$ and the last condition is the consequence of the fact that $w^{-1}\equiv 0$ at $t=0$ (since the circle acts trivially on $X$).
  \par
 Our first result, which is a singular Cauchy-Kovalevskaya theorem, says that we can indeed solve this system locally, if the initial data $h_{ij}$ is real-analytic.
\begin{theorem} Let $h_{ij}$, $i,j=1,\dots, n$, be real-analytic functions on an open subset $U$ of ${\Bbb C}^n$. Then there exits a unique solution of the system \eqref{system} on an open neighbourhood of $U$ in $U\times [0,+\infty)$.\label{exists}\end{theorem}
\begin{remark} Observe that, if the solution does give the metric on a neighbourhood of $X$ in $L$, then the initial data $h_{ij}$ is real-analytic, since the metric $g$ is Ricci-flat.\end{remark}
\begin{pf} We treat ${\Bbb C}^n$ as a real subspace $V$ of ${\Bbb C}^{2n}$. Since $h_{ij}$ are real analytic, they extend to holomorphic functions on a neighbourhood of $V$. Therefore we can treat problem \eqref{system} as purely holomorphic, i.e. $x_i,y_j$ are complex coordinates. First of all, it is easy to see that \eqref{system} has a unique formal solution, i.e. a power series in $t$. Thus we only have to show that this series is convergent. Let us write 
 $e^u=te^v$. Then we can rewrite the problem \eqref{system} as 
 \begin{equation} \begin{cases} tv_t=-1+ce^{-v}\det[g_{ij}]\\ v_{x_i{x}_j}+v_{y_i{y}_j}+c(g_{ij})_t=0\end{cases},\end{equation}
 with the initial conditions $(g_{ij})_{|_{t=0}}=h_{ij}$, $\left(e^v\right)_{|_{t=0}} = c\det h$.
 A theorem showing convergence of a formal solution to this system is proved in the appendix. This theorem is a slight generalization of a theorem of G\'erard and Tahara \cite{G-T}. It is applied to functions $\tilde{g}_{ij}=g_{ij}-h_{ij}$ and $\tilde{v}=v-v_0$, where $v_0=v_{|_{t=0}}$.\end{pf}

Having solved the Cauchy problem \eqref{system} we ask whether the solution gives us a smooth metric \eqref{metric} on a neighbourhood of $X$ in some line bundle $L$. First of all, we have
\begin{lemma} Suppose that we have a local solution of the Cauchy problem \eqref{system} (with $c\neq 0$). Then the metric \eqref{metric} extends smoothly to the hypersurface $t=0$ (which is the fixed-point set of the circle action) if and only if $c=1$.\end{lemma}
\begin{pf} 
Since $\det g$ is finite and non-zero at $t=0$, equation \eqref{toda2} implies that $e^u=t(a+bt+\ldots)$ near $t=0$ with $a\neq 0$. Therefore $u_t=\frac{1}{t}+O(1)$ near $t=0$. Now $w=c^{-1}u_t$. This implies immediately that $c$ must be positive. Furthermore, the fibrewise metric is 
$$wdt^2+w^{-1}\phi^2=\left(\frac{1}{ct}+O(1)\right)dt^2+ \left(\frac{1}{ct}+O(1)\right)^{-1}\phi^2.$$
If we introduce a new coordinate $r$ so that $t=r^2$, we see that this metric extends smoothly to the origin if and only if $c=1$. Now, the formula \eqref{curvature} shows that the connection $1$-form $\phi$ and hence the metric \eqref{metric} extends to the hypersurface $t=0$. 
\end{pf}
\medskip

Thus it remains to show that the curvature form \eqref{curvature} belongs to $2\pi{\Bbb Z}$.  We observe that 
the first equation in \eqref{system} (with $c=1$) says that
\begin{equation}\frac{d}{dt}\omega=-i\partial \bar{\partial}u,\label{o1} \end{equation}
where $\omega=\omega(t)$ is the K\"ahler form of the metric $g$ at time $t$ on $X$. On the other hand $u=\ln\left(w^{-1} \det [g_{ij}]\right)$.  Since $g$ is K\"ahler \eqref{o1} can be written as
\begin{equation} \frac{d}{dt}\omega=\rho\left(g\right)+ i\partial\bar{\partial}\ln w \label{ricci} \end{equation}
where $\rho$ denotes the Ricci form of a K\"ahler metric. Now \eqref{curvature} shows that the curvature form of $\phi$ is indeed in $2\pi {\Bbb Z}$, and in fact represents $-c_1(X)$. Thus Theorem 1 is proved.

\section{Examples}

We wish now to give explicit examples of Ricci-flat K\"ahler metrics on canonical bundles. Given a K\"ahler manifold $(X,h,I)$ we seek a time dependent metric $g=g(t)$ and a function $w$ on $X\times I$
which satisfy the equations \eqref{ricci}, \eqref{cw} and \eqref{u}. The last two give us
\begin{equation} w^{-1}=\frac{\int_0^t \det g}{\det g}. \label{w}\end{equation}
Substituting into \eqref{ricci} we obtain
\begin{equation} \frac{d}{dt}\omega=-i\partial\bar{\partial} \int_0^t \det g\label{nricci} \end{equation}

 \par
The first example deals with manifolds with constant principal Ricci curvatures, i.e. constant eigenvalues of the Ricci curvature. This class of manifolds includes both  K\"ahler-Einstein manifolds and homogeneous manifolds. The following result is a particular case of a theorem of Hwang and Singer \cite{HS}. 

\begin{theorem} Let $X^{2n}$ be a K\"ahler manifold with K\"ahler form $\Phi$ such that the eigenvalues of the Ricci curvature are constant. Then the solution to \eqref{nricci} is given by
\begin{equation} \omega=\Phi+t\rho(\Phi).\label{omega-hom}\end{equation}
The function $w^{-1}$ is given by
\begin{equation} w^{-1}(t)= \frac{\int_0^t P(t)dt}{P(t)}, \label{K-hom}\end{equation}
where $P(t)$ is defined as $P(t)=\left(\Phi+t\rho(\Phi)\right)^n/\Phi^n$ (and so it depends only on the eigenvalues of the Ricci curvature). 
\par
In particular, if all the eigenvalues of the Ricci curvature are nonnegative, then the resulting Ricci-flat metric on $K_X$ is complete.\end{theorem}

\begin{pf} It is sufficient to observe that $\omega^n=P(t)\Phi^n$, so $\rho(\omega)=\rho(\Phi)$. Now it is clear that $\omega$  satisfies \eqref{nricci}. \end{pf}

We remark that Apostolov, Armstrong and Draghici \cite{AAD}  recently found examples of irreducible non-homogeneous K\"ahler manifolds with constant principal Ricci curvatures. 

As an aside, let us give an application to the geometry of K\"ahler quotients. We recall that Futaki \cite{Futaki} has shown that if $M$ is a compact K\"ahler-Einstein manifold with positive scalar curvature and a Hamiltonian Killing vector field whose length is constant on the level sets of the moment map, then the K\"ahler quotient by the resulting circle action is also K\"ahler-Einstein. A simple example of ${\Bbb C}\times {\Bbb C}^2$ with the diagonal circle action on the second factor (and trivial on the first) shows that Futaki's result does not hold for Ricci-flat manifolds. Nevertheless we have a weaker conclusion.

\begin{proposition} Let $M$ be a complete Ricci-flat K\"ahler manifold with an isometric and Hamiltonian circle action such that the length of the Killing vector field is constant on the level sets of the moment map. Moreover, assume that the moment map is bounded from  below. Then the K\"ahler quotient of $M$ by $S^1$ has constant principal Ricci curvatures.\end{proposition}

\begin{pf} Let $t$ be the moment map and $X=t^{-1}(a)/S^1$ a particular K\"ahler quotient of $M$. Since $a$ is a regular value of $t$, the K\"ahler quotients for nearby level sets are isomorohic tp $X$ (as a complex manifold). Thus we have a family $g(t)$ of K\"ahler metrics on $X$. The assumption and the equation \eqref{new} show that $g(t)$ is linear in $t$. Since $t$ is bounded from below, Proposition \ref{Laplacian} implies that the constant $c$ of Theorem \ref{P-P} is non-zero. Now the equations \eqref{cw} and \eqref{u} show that $\det g(t)=f(t)\det g(a)$ for $t$ near $a$. Therefore the Ricci form $\rho(g)$ is constant in $t$ and the equation  \eqref{PP} implies that $\rho(g)$ is equal to $\omega_t$, where  $\omega$ is the K\"ahler form of $g$.  The conclusion follows now, since we already know that  $\omega(t)^n=f(t)\omega(a)^n$. \end{pf}
\medskip 
 
The next examples involve surfaces of revolution. Let $\Sigma$ denote either ${\Bbb C}$ or ${\Bbb C}P^1$ with a metric of constant curvature. We define a surface $\Sigma_a$ as the K\"ahler quotient of $\Sigma \times {\Bbb C}$ by the action of ${\Bbb R}$ defined as
$$ r\times (x,z)=(e^{2\pi iar}\cdot x, z+r).$$
This is a surface of revolution and the Ricci-flat K\"ahler metric on $K_{\Sigma_a}=T^\ast \Sigma_a$ is complete, since it can be obtained as a hyperk\"ahler quotient of $T^\ast \Sigma\times {\Bbb H}$ by ${\Bbb R}$. Now, the main result of \cite{Biel} shows that these are all such surfaces of revolution:

\begin{proposition} Let $(X,h)$ be a surface of revolution such that the Ricci-flat K\"ahler metric defined in Theorem \ref{main} is complete. Then $(X,h)$ is isometric to one of the surfaces $\Sigma_a$ defined above. \end{proposition}

\bigskip

\appendix
\section{A Cauchy-Kovalevskaya theorem for a class of singular PDE's}
In this section we shall prove a result about convergence of formal solutions to certain singular partial differential equations. This is a generalization of a theorem of G\'erard and Tahara \cite{G-T} to a class of singular systems of PDE's and it uses their method of proof. The result of G\'erard and Tahara applies to first order nonlinear PDE's of the form:
$$ g\left(t,x,v,t\frac{\partial v}{\partial t}, \frac{\partial v}{\partial x_1},\dots,\frac{\partial v}{\partial x_n}\right)=0$$
where $g$ is a holomorphic function in some polydisc. As observed in \cite{G-T}, the convergence of formal solutions is no longer true if we allow $g$ to depend on second derivatives of $v$. 
We shall now show that the theorem remains valid for systems of PDE's of the following form:
\begin{equation}\begin{cases}g\left(t,x,v,t\frac{\partial v}{\partial t}, \frac{\partial v}{\partial x_1},\dots,\frac{\partial v}{\partial x_n}, w_1,\dots,w_N\right)=0 \\ \frac{\partial w_i}{\partial t}=L_i(x)(v)+a_i(t,x), \quad\quad i=1,\dots,N \end{cases}\label{E0}\end{equation}
where $L_i(x)$ are linear differential operators of order at most $2$. We remark that one can allow the dependence of $L_i$ on $t$, but this further complicates the already complicated notation.
\par
To guarantee the existence of formal solutions we shall assume that the first equation can be written as:
$$\left(t\frac{\partial}{\partial t} -\rho(x)\right)v= tb(x)+G(x)\left(t,v,t\frac{\partial v}{\partial t}, \frac{\partial v}{\partial x_1},\dots,\frac{\partial v}{\partial x_n}, w_1,\dots,w_N\right)$$
where $\rho(x)$ and $b(x)$ are holomorphic functions defined in a polydisc $D$ centered at the origin of ${\Bbb C}^n$, and
$$G(x)(t,Z,V,X_i, Y_{j})_{i\leq n, j\leq N}=\sum_{p+r+s+|\alpha|+2|\beta|\geq 2} a_{p,r,s,\alpha,\beta}(x)  t^pZ^q V^s X_1^{\alpha_1}\dots X_n^
{\alpha_n}Y_{1}^{\beta_{1}} \dots Y_{N}^{\beta_{N}}.$$
The coefficients $a_{p,q,s,\alpha,\beta}(x)$ are holomorphic in $D$ and 
$$|a_{p,r,s,\alpha,\beta}(x)|\leq A_{p,q,s,\alpha,\beta}.$$
Moreover the power series
$$ \sum A_{p,q,s,\alpha,\beta}t^pZ^qU^sX^\alpha Y^\beta$$
is convergent near the origin.
\par
We seek a holomorphic solution $(v,w_i)$ to the above system satisfying
$$ v(0,x)=w_i(0,x)\equiv 0, \quad i=1,\dots,N.$$
A formal solution is a power series solution of the form
$$\sum_{m\geq 1} v_m(x)t^m $$
whose coefficients are holomorphic in $D$.

The particular form of the system \eqref{E0} allows us to rewrite it as a single differential-integral equation:
\begin{equation} \left(t\frac{\partial}{\partial t} -\rho(x)\right)v= tb(x)+G(x)\left(t,v,t\frac{\partial v}{\partial t}, \frac{\partial v}{\partial x_1},\dots,\frac{\partial v}{\partial x_n},\int_0^tL_1(v),\dots,\int_0^tL_N(v)\right).\label{E}\end{equation}
Here we regrouped the terms in the power expansion of $G$, so that $G$ does not depend on $\int_0^t a_i(,t,x)$.
\begin{theorem} Each formal solution of \eqref{E} is convergent. If $\rho(0)\not\in{\Bbb N}^\ast$, then there exists a
unique formal solution satisfying $v(0,x)\equiv 0$. \label{C-K}\end{theorem}
\begin{pf}
If $\rho(0)\not\in{\Bbb N}^\ast$, then \eqref{E} has a unique formal solution of the form
\begin{equation}\sum_{m\geq 1} v_m(x)t^m.\label{um}\end{equation}
Moreover, $v_m(x)$ is determined recursively by the following formula:
$$v_1(x)=\frac{b(x)}{1-\rho(x)},$$
and for $m\geq 2$
\begin{multline}v_m(x)=\frac{1}{m-\rho(x)}f_m\Bigl(v_1,2v_2,\dots, (m-1)v_{m-1},v_1,\dots,v_{m-1},\partial_1 v_1, 
\dots, \partial_n v_1,\dots,\\\partial_1 v_{m-1},\dots,\partial_n v_{m-1}, \frac{1}{2}L_1(v_1)
,\dots,\frac{1}{2}L_N(v_1),\frac{1}{m}L_1(v_{m-1}),\dots, \frac{1}{m}L_N(v_{m-1});\\ 
\left\{a_{p,q,s,\alpha,\beta}\right\}_{p+q+s+|\alpha|+|\beta|
\leq m}\Bigr).\label{fm}\end{multline}

We shall show that this solution is convergent. Let $D_a$ denote the polydisc of diameter $2a$. By taking $R$ sufficiently
small (in particular, $R<1$), we can assume that all the $v_m(x)$ are holomorphic in $D_R$ and we have;
$$|v_1(x)|\leq A, \enskip |\partial_i v_1 (x)|\leq A,\enskip i=1,\dots, n\enskip |L_j(v_1) (x)|\leq A,\enskip j=1,\dots,N;$$
$$|m-\rho(x)|\geq \sigma m,\enskip m=1,2,3,\dots.$$
Moreover, let $M$ be a constant such that the coefficients of
$$L_i=\sum c^i_{kl}(x)\frac{\partial^2}{\partial x_k\partial x_l}+\sum d^i_k(x)\frac{\partial }{\partial x_k}+e(x)$$
satisfy 
$$ \sum \left|c^i_{kl}(x)\right|+\sum \left|d^i_k(x)\right|+|e(x)|\leq M.$$

Now we consider the analytic equation:
\begin{multline*}\sigma Y=\sigma At+\\\frac{1}{(R-r)^2}\sum_{p+q+s+|\alpha|+2|\beta|\geq
2}\frac{A_{p,q,s,\alpha,\beta}}{(R-r)^{p+q+s+|\alpha|+2|\beta|-2}} t^p Y^{q+s}(2eY)^{\bigl(\sum\alpha_i\bigr)} 
(4e^2MYt)^{\bigl(\sum\beta_{i}\bigr)}.\end{multline*}
Here $e$ is the smallest real number such that $e^{\pi \sqrt{-1}}=-1$.

By the implicit function theorem, this equation has a unique analytic solution of the form 
$$Y=\sum_{m\geq 1}Y_m(r)t^m,$$
determined by the following recursive formula
$$Y_1=A$$
and, for $m\geq 2$,
\begin{multline} \sigma Y_m=\frac{1}{(R-r)^2}F_m\Bigl(Y_1,\dots,Y_{m-1},2eY_1,\dots, 2eY_{m-1},4e^2MY_1,\dots,
4e^2MY_{m-1};\\ 
\left.\left\{\frac{A_{p,q,s,\alpha,\beta}}{(R-r)^{p+q+s+|\alpha|+2|\beta|-2}}\right\}_{p+q+s+|\alpha|+|\beta|\leq m}\right).
\label{Fm}\end{multline}

Moreover, by induction on $m$, we see that $Y_m(r)$ is expressed in the form 
\begin{equation} Y_m(r)=\frac{C_m}{(R-r)^{2m-2}},\enskip m=1,2,\dots,\label{Ym}\end{equation}
with constants $C_1=A$ and $C_m\geq 0$ (for $m\geq 2$).
\par
We shall show that the power series for $Y$ is a majorant power series for the formal solution \eqref{um}. To do so, it is
sufficient to prove the following inequalities for all $m$:
\begin{equation} |v_m(x)| \leq m|v_m(x)|\leq Y_m(r)\enskip \text{on $D_r$ for $0<r<R$};\label{first}\end{equation}
\begin{equation} \left|\partial_i v_m(x)\right|\leq 2eY_m(r)\enskip \text{on $D_r$ for $0<r<R$}, i=1,\dots,n;\label{second}
\end{equation}
\begin{equation} \left|L_k\left(v_m\right)(x)\right|\leq 4e^2(m+1)MY_m(r)\enskip \text{on $D_r$ for $0<r<R$}, \enskip k=1,\dots,N.
\label{third}
\end{equation}

The case $m=1$ is clear from the definition of $A$. We proceed by induction. We replace all the terms in \eqref{fm} by
their absolute values. Then we use the inductive assumption and also replace $|a_{p,q,s,\alpha,\beta}|$ by 
$\frac{A_{p,q,s,\alpha,\beta}}{(R-r)^{p+q+s+|\alpha|+2|\beta|-2}}$ (this is a majorant, as $R<1$). This has also the 
effect of replacing $f_m$ by $F_m$, and, from \eqref{Fm}, it gives:
\begin{equation} |v_m(x)|\leq \frac{1}{m}(R-r)^2 Y_m(r)\label{First}\end{equation}
which proves \eqref{first} (cf. \cite{G-T}, p. 985). Since $Y_m(r)$ has the form \eqref{Ym}, the above inequality can be 
written as:
$$|v_m(r)|\leq \frac{1}{m}\frac{C_m}{(R-r)^{2m-4}}.$$
Now the following lemma  proves \eqref{second} and \eqref{third}.
\begin{lemma} If a function $v(x)$ holomorphic in $D_R$ satisfies 
$$|v(x)|\leq\frac{C}{(R-r)^p}\enskip \text{on $D_r$ for $0<r<R$},$$
then
$$\left|\partial_i v(x)\right|\leq\frac{Ce(p+1)}{(R-r)^{p+1}}\enskip \text{on $D_r$ for $0<r<R$, $i=1,\dots,n$}.$$
\end{lemma}
 For the proof see \cite{Hor}, Lemma 5.1.3.
 \medskip

We now assume that $\rho(0)=k\in {\Bbb N}^\ast$. We can modify
$\left\{A_{p,q,s,\alpha,\beta}\right\}_{p+q+s+|\alpha|+|\beta|\leq k}$ so that $v_k(x)$ satisfies
\eqref{first}-\eqref{third} and then apply the previous proof.

\end{pf}

\begin{ack} The author thanks V. Apostolov, O. Biquard, D. Calderbank, M. Ro\v{c}ek and M. Singer for useful discussions and comments.\end{ack}


\begin{thebibliography}{99}

\bibitem{AAD} 
{V. Apostolov, J. Armstrong and T. Draghici}, `Local rigidity of certain classes of almost Kahler 4-manifolds', math/9911197.

\bibitem{Biel}
{R. Bielawski}, `Complete hyper-Kähler $4n$-manifolds with a local tri-Hamiltonian $\bold R\sp
n$-action', {\it Math. Ann.} 314 (1999), 505--528.

\bibitem{Cal}
{E. Calabi}, `M\'etriques k\"ahl\'eriennes et fibr\'es holomorphes',  
{\it Ann. Sci. \'Ecole Norm. Sup. (4)} 12 (1979), 269--294.

 
\bibitem{Feix}
{B. Feix},  University of Cambridge Ph.D. Thesis, Cambridge (2000).

\bibitem{Futaki}
{A. Futaki}, `The Ricci curvature of symplectic quotients of Fano manifolds', 
{\it Tôhoku Math. J. (2)} 39 (1987), 329--339. 

\bibitem{G-T} 
{R. G\'erard \and H. Tahara},  `Holomorphic and singular solutions of nonlinear singular first order partial
differential equations', {\it Publ. RIMS, Kyoto Univ.} 26 (1990), 979--1000.

\bibitem{Hor}
{L. H\"ormander}, {\em Linear partial differential operators}, Springer  (1963).

\bibitem{HS} 
{A. Hwang and M. Singer}, ` A momentum construction for circle-invariant Kahler metrics', math/9811024.

\bibitem{Kha}
{D. Kaledin}, `Hyperkaehler structures on total spaces of holomorphic cotangent bundles', alg-geom/9710026.

\bibitem{LeB}
{C. LeBrun}, `Explicit self-dual metrics on ${\Bbb C}P^2\#\cdots\#{\Bbb C}P^2$',  {\it J. Differential Geom.}
34 (1991), 223--253.


\bibitem{PP}
{H. Pedersen \and Y.S. Poon}, `Hamiltonian constructions of Kähler-Einstein metrics and Kähler metrics of
constant scalar curvature', {\it Comm. Math. Phys.} 136 (1991), 309--326.




\end{thebibliography}
\end{document}